\overfullrule=0pt
\centerline {\bf Miscellaneous applications of certain minimax theorems II}\par
\bigskip
\bigskip
\centerline {BIAGIO RICCERI}\par
\bigskip
\bigskip
\centerline {\it Dedicated to Professor Hoang Tuy with my greatest esteem}\par
\bigskip
{\bf Abstract.} In this paper, we present new applications of our general minimax theorems. In particular, one of them concerns the
multiplicity of global minima for the integral functional of the Calculus of Variations.\par
\bigskip
{\bf Keywords:} minimax theorem; connectedness; global minimum; multiplicity; integral functional; Neumann problem.\par
\bigskip
{\bf Mathematics Subject Classification (2010):} 49J35; 49K35; 49J45; 49K27; 35J92; 90C47.
\bigskip
\bigskip
\bigskip
\bigskip
{\bf 1. Statements of the main results}\par
\bigskip
This is the second ring of a chain of papers (started with [14]) which is devoted to consequences and applications of certain general
minimax theorems that we have established in the past years ([5]-[15]).\par
\smallskip
The motivation for such papers is just to show the great flexibility and usefulness of those theorems.\par
\smallskip
The two main results that we want to prove in the present paper are Theorems 1.1 and 1.2 below.\par
\medskip
A real-valued function $f$ on a topological space is said to be inf-connected (resp. sup-connected)
if $f^{-1}(]-\infty,r[)$ (resp. $f^{-1}(]r,+\infty[)$) is connected for all $r\in {\bf R}$.\par
\medskip
THEOREM 1.1. - {\it Let $X, Y$ be two real Banach spaces, $\Phi:X\to Y$ a surjective continuous linear operator, $\Psi:X\to Y$ a non-constant Lipschitzian operator with Lipschitz constant equal to $L$, $\varphi:Y\to {\bf R}$ a non-constant, continuous, concave and inf-connected functional, $[a,b]$ a closed sub-interval of $[-1,1]$.\par
Then, for every continuous and concave function $\gamma:[a,b]\to {\bf R}$, one has
$$\max\left \{
\inf_{x\in X}\varphi\left (\Phi(x)+{{a}\over {\alpha_ {\Phi}L}}\Psi(x)\right )+\gamma(a),
\inf_{x\in X}\varphi\left (\Phi(x)+{{b}\over {\alpha_ {\Phi}L}}\Psi(x)\right )+\gamma(b)\right \}$$
$$=\inf_{x\in X}\sup_{\lambda\in [a,b]}\left (\varphi\left (\Phi(x)+{{\lambda}\over {\alpha_ {\Phi}L}}\Psi(x)\right )+\gamma(\lambda)\right )\ ,$$
where
$$\alpha_{\Phi}=\sup_{\|y\|_Y\leq 1}\hbox {\rm dist}(0,\Phi^{-1}(y))\ .$$} \par
\medskip
 Let $\Omega\subset {\bf R}^n$ be a bounded domain with smooth boundary and let $p>1$. On the Sobolev space $W^{1,p}(\Omega)$, we
consider the norm
$$\|u\|=\left ( \int_{\Omega}|\nabla u(x)|^p dx+\int_{\Omega}|u(x)|^p dx\right ) ^{1\over p}\ .$$
If $n\geq p$, we denote by ${\cal A}$ the class of all
continuous functions $\psi:{\bf R}\to {\bf R}$ such that
$$\sup_{\xi\in {\bf R}}{{|\psi(\xi)|}\over
{1+|\xi|^q}}<+\infty\ ,$$
where  $0<q< {{pn}\over {n-p}}$ if $p<n$ and $0<q<+\infty$ if
$p=n$. While, when $n<p$, ${\cal A}$ stands for the class
of all continuous functions $\psi:{\bf R}\to {\bf R}$.\par
\smallskip
Recall that a function $\varphi:\Omega\times {\bf R}^{m}\to {\bf R}$
 is said to be a normal integrand ([16]) if it is
${\cal L}(\Omega)\otimes {\cal B}({\bf R}^m)$-measurable and
$\varphi(x,\cdot)$ is lower semicontinuous for a.e. $x\in \Omega$.
Here ${\cal L}(\Omega)$ and ${\cal B}({\bf R}^m)$ denote the
Lebesgue and the Borel $\sigma$-algebras of subsets of $\Omega$ and
${\bf R}^m$, respectively. \par
\smallskip
Recall that if $\varphi$ is a normal integrand then, for each measurable function $u:\Omega\to
{\bf R}^m$, the composite function $x\to \varphi(x,u(x))$ is measurable
([16]).\par
\smallskip
A real-valued function $f$ on a convex set is said to be quasi-convex (resp. quasi-concave) if $f^{-1}(]-\infty,r[)$ (resp.
$f^{-1}(]r,+\infty[)$) is convex for all $r\in {\bf R}$.\par
\medskip
THEOREM 1.2. - {\it  Let
$\varphi:\Omega\times {\bf R}\times {\bf R}^n\to {\bf R}$ be a normal integrand such that $\varphi(x,\xi,\cdot)$
is convex for all $(x,\xi)\in\Omega\times {\bf R}$ and let $\psi\in {\cal A}$ be a strictly monotone function. Assume that:\par
\noindent
$(i)$\hskip 5pt there are $c, d>0$ such that
$$c|\eta|^p-d\leq \varphi(x,\xi,\eta)$$
for all $(x,\xi,\eta)\in \Omega\times{\bf R}\times{\bf R}^n$ and
$$\lim_{|\xi|\to +\infty}{{\inf_{(x,\eta)\in \Omega\times{\bf R}^n}\varphi(x,\xi,\eta)}\over {|\psi(\xi)|+1}}=+\infty\ ;$$
\noindent
$(ii)$\hskip 5pt for each $\xi\in {\bf R}$, the function $\varphi(\cdot,\xi,0)$ lies in $L^1(\Omega)$ and the function
$\int_{\Omega}\varphi(x,\cdot,0)dx$ is not quasi-convex.\par
Then, for every sequentially weakly closed set $V\subseteq W^{1,p}(\Omega)$, containing the constants, and for every convex set
$Y\subseteq L^{\infty}(\Omega)$, dense in $L^{\infty}(\Omega)$, there exists $\alpha\in Y$
such that the restriction to $V$ of the functional 
$$u\to \int_{\Omega}\varphi(x,u(x),\nabla u(x))dx+\int_{\Omega}\alpha(x)\psi(u(x))dx$$
 has at least two global minima. The same property holds also with $Y=C^{\infty}_0(\Omega)$.}\par
\bigskip
{\bf 2. Tools for proving the main results}\par
\bigskip
In this section, for reader's convenience, we collect the tools that we will use to prove Theorems 1.1 and 1.2. For the basic notions
on multifunctions, we refer to [1]. The results without any reference are new.
\medskip
THEOREM 2.A ([10], Theorem 5.7). - {\it Let $X$ be a topological space, $I$ a compact real interval and $f:X\times I\to {\bf R}$ a function
which is lower semicontinuous in $X$, and quasi-concave and upper semicontinous in $I$. Moreover, assume that the set
$$\{\lambda\in I : f(\cdot,\lambda)\hskip 3pt \hbox {\it is\hskip 3pt inf-connected\hskip 3pt in}\hskip 3pt X\}$$
is dense in $I$.\par
Then, one has 
$$\sup_I\inf_Xf=\inf_X\sup_If\ .$$}\par
\medskip
A real-valued function $f$ on a topological space is said to be inf-compact if $f^{-1}(]-\infty,r])$ is compact for all $r\in {\bf R}$.\par
\medskip
THEOREM 2.B ([15], Theorem 1.2). - {\it Let $X$ be a topological space, $E$ a real vector
space, $Y\subseteq E$ a non-empty convex set and $f:X\times Y\to {\bf R}$ a function which is lower semicontinuous and inf-compact in $X$, 
and concave in $Y$. Moreover, assume that
$$\sup_Y\inf_Xf<\inf_X\sup_Yf\ .$$
Then, there exists $\hat y\in Y$ such that the function $f(\cdot,\hat y)$ has at least two global
minima.}
\medskip
PROPOSITION 2.A ([10], Proposition 5.6). - {\it Let $X, Y$ be two topological spaces,
$F:X\to 2^Y$ a lower semicontinuous multifunction with non-empty values and $A\subseteq X$ a connected set. Assume that the set
$$\{x\in A : F(x)\hskip 3pt \hbox {\it is\hskip 3pt connected}\}$$
is dense in $A$.\par
Then, the set $F(A)$ is connected.}\par
\medskip
THEOREM 2.C ([4], Th\'eor\`eme 2). - {\it  Let $X, Y$ be two real Banach spaces, $\Phi:X\to Y$ a surjective continuous linear operator, $\Psi:X\to Y$ a non-constant Lipschitzian operator with Lipschitz constant equal to $L$. Assume that $\alpha_{\Phi}L<1$.\par
Then, the multifunction $y\to (\Phi+\Psi)^{-1}(y)$ is Lipschitzian in $Y$ and its values are absolute retracts.}\par
\medskip
We now are in a position to prove the following result from which we will draw Theorem 1.1:
\medskip
THEOREM 2.1. - {\it Let $X$ be a topological space, $E$ a real topological vector space, $Y$ a convex subset of $E$, $\varphi:Y\to {\bf R}$ a continuous, concave and inf-connected functional, $I$ a compact real interval, 
$f, g:X\to E$ two continuous functions such that $f(x)+\lambda g(x)\in Y$ for all $x\in X$,
$\lambda\in I$.  Moreover,
assume that there exists a set $D\subseteq I$, dense in $I$, with the following property:  for each $\lambda\in D$, the function
$f+\lambda g$ is onto $Y$ and open with respect to the relative topology of $Y$,  and there exists a set $S_{\lambda}\subseteq Y$, dense in $Y$, such that
the set $(f+\lambda g)^{-1}(y)$ is connected for each $y\in S_{\lambda}$.\par
Then, for every continuous and concave function $\gamma:I\to {\bf R}$, one has
$$\inf_{x\in X}\sup_{\lambda\in I}(\varphi(f(x)+\lambda g(x))+\gamma(\lambda))=\sup_{\lambda\in I}\inf_{x\in X}(\varphi(f(x)+\lambda g(x))+\gamma(\lambda))\ .$$} 
\smallskip
PROOF. Consider the function $\psi:X\times I\to {\bf R}$ defined by
$$\psi(x,\lambda)=\varphi(f(x)+\lambda g(x))+\gamma(\lambda)$$
for all $(x,\lambda)\in X\times I$. Clearly, for each $x\in X$, the function $\psi(x,\cdot)$ is concave and continuous in $I$. Now, fix $\lambda\in D$.
Let $r\in {\bf R}$ be such that $\{x\in X : \psi(x,\lambda)<r\}\neq\emptyset$.
 Clearly, we have
$$\{x\in X : \psi(x,\lambda)<r\}=(f+\lambda g)^{-1}(\varphi^{-1}(]-\infty,r-\gamma(\lambda)[))\ .$$
Now, observe that $\varphi^{-1}(]-\infty,r-\gamma(\lambda)[)$ is open in $Y$ and connected since $\varphi$ is continuous and inf-connected. But, since $\lambda\in D$, the multifunction $y\to (f+\lambda g)^{-1}(y)$ is non-empty valued and lower semicontinuous in $Y$. Since $S_{\lambda}\cap \varphi^{-1}(]-\infty,r-\gamma(\lambda)[)$ is dense in
$\varphi^{-1}(]-\infty,r-\gamma(\lambda)[)$, thanks to Proposition 2.A,
we conclude that the set
 $(f+\lambda g)^{-1}(\varphi^{-1}(]-\infty,r-\gamma(\lambda)[))$ is connected. Clearly, $\psi(\cdot,\lambda)$ is continuous in $X$ for all $\lambda\in I$.
Now, the conclusion follows directly from Theorem 2.A.\hfill $\bigtriangleup$\par
\medskip
The following two results will be used jointly with Theorem 2.B to prove Theorem 1.2.\par
\medskip
PROPOSITION 2.1. - {\it Let $\Omega\subset {\bf R}^n$ be a bounded domain with smooth boundary, let $p>1$ and let
$\varphi:\Omega\times {\bf R}\times {\bf R}^n\to {\bf R}$ be normal integrand such that, for some $c, d>0$, one has
$$c|\eta|^p-d\leq \varphi(x,\xi,\eta)$$
for all $(x,\xi,\eta)\in\Omega\times{\bf R}\times{\bf R}^n$ and
$$\lim_{|\xi|\to +\infty}\inf_{(x,\eta)\in \Omega\times{\bf R}^n}\varphi(x,\xi,\eta)=+\infty\ .$$
Then, in $W^{1,p}(\Omega)$, one has
$$\lim_{\|u\|\to +\infty}\int_{\Omega}\varphi(x,u(x),\nabla u(x))dx=+\infty\ .$$}\par
\smallskip
PROOF.  Clearly, for each $u\in W^{1,p}(\Omega)$, 
the integral $\int_{\Omega}\varphi(x,u(x),\nabla u(x))dx$ exists
and belongs to $]-\infty,+\infty]$. Fix a sequence $\{u_n\}$
in $W^{1,p}(\Omega)$ such that $\lim_{n\to \infty}\|u_n\|=+\infty$.
We have to show, up to a sub-sequence, that
$$\lim_{n\to \infty}\int_{\Omega}\varphi(x,u_n(x),\nabla u_n(x))dx =+\infty\ .$$
If the sequence $\left \{ \int_{\Omega}|\nabla u_n(x)|^p dx\right \}$
is unbounded, this clearly holds, due to the assumed growth of $\varphi$.
 So, assume that the sequence $\left \{ \int_{\Omega}|\nabla u_n(x)|^p dx\right \}$
is bounded.
 Then, by the Poincar\'e-Wirtinger
inequality, there exists a constant $\gamma>0$ such that
$$\int_{\Omega}|u_n(x)-a_n|^p dx\leq \gamma$$
for all $n\in {\bf N}$, where
$$a_n={{\int_{\Omega}u_n(x)dx}\over {m(\Omega)}}\ ,$$
$m(\Omega)$ being the Lebesgue measure of $\Omega$.
Since
$$\lim_{n\to \infty}\int_{\Omega}|u_n(x)|^p dx=+\infty$$
we clearly have
$$\lim_{n\to \infty}|a_n|=+\infty\ .$$
Fix any $M>0$. Then, there exists $\delta>0$ such that
$$\varphi(x,\xi,\eta)\geq M$$
for all $(x,\xi,\eta)\in \Omega\times{\bf R}\times{\bf R}^n$ with $|\xi|\geq\delta$.
We now show that
$$\lim_{n\to \infty}m(A_n)=m(\Omega)$$
where
$$A_n=\{x\in \Omega : |u_n(x)|\geq \delta\}\ .$$
Arguing by contradiction, assume that
$$\liminf_{n\to \infty}m(A_n)<m(\Omega)\ .$$
Fix $\rho$ satisfying
$$\liminf_{n\to \infty}m(A_n)<\rho<m(\Omega)\ .$$
Now, fix $\theta>\gamma$ and $n\in {\bf N}$ so that, at the same time, one has
$$|a_n|>\left ( {{\theta}\over {m(\Omega)-\rho}}\right ) ^{1\over p}+\delta$$
as well as
$$m(A_n)<\rho\ .$$
Then, one has
$$\gamma\geq \int_{\Omega}|u_n(x)-a_n|^p dx\geq \int_{\Omega\setminus A_n}
|u_n(x)-a_n|^p dx> (|a_n|-\delta)^p (m(\Omega)-\rho)>\theta\ ,$$
an absurd. Now, for each $n\in {\bf N}$, we have
$$\int_{\Omega}\varphi(x,u_n(x),\nabla u_n(x))dx=\int_{A_n}\varphi(x,u_n(x),\nabla u_n(x))dx+
\int_{\Omega\setminus A_n}\varphi(x,u_n(x),\nabla u_n(x))dx\geq$$
$$ Mm(A_n)-m(\Omega\setminus A_n)d\ .$$
Therefore
$$\liminf_{n\to \infty}\int_{\Omega}\varphi(x,u_n(x),\nabla u_n(x))dx\geq Mm(\Omega)\ .$$
Since $M$ is arbitrary, the sequence $\{\int_{\Omega}\varphi(x,u_n(x),\nabla u_n(x))dx\}$ diverges and the proof
is complete.\hfill $\bigtriangleup$
\par
\medskip
PROPOSITION 2.2. - {\it Let $X, Y$ be two non-empty sets and $f:X\to {\bf R}$, $g:X\times Y\to {\bf R}$ two given
functions. Assume that there are two sets $A, B\subset X$ such that:\par
\noindent
$(a)$\hskip 5pt $\sup_Af<\inf_Bf$\ ;\par
\noindent
$(b)$\hskip 5pt $\sup_{y\in Y}\inf_{x\in A}g(x,y)\leq 0$\ ;\par
\noindent
$(c)$\hskip 5pt $\inf_{x\in B}\sup_{y\in Y}g(x,y)\geq 0$\ ;\par
\noindent
$(d)$\hskip 5pt $\inf_{x\in X\setminus B}\sup_{y\in Y}g(x,y)=+\infty$\ .\par
Then, one has
$$\sup_{y\in Y}\inf_{x\in X}(f(x)+g(x,y))<\inf_{x\in X}\sup_{y\in Y}(f(x)+g(x,y))\ .$$}\par
\smallskip
PROOF. Fix $y\in Y$ and $\epsilon\in ]0,\inf_Bf-\sup_Af[$ as well. Since $\inf_{x\in A}g(x,y)\leq 0$, there
is $\tilde x\in A$ such that $g(\tilde x,y)<\epsilon$. Hence, we have
$$\inf_{x\in X}(f(x)+g(x,y))\leq f(\tilde x)+g(\tilde x,y)<\sup_Af+\epsilon\ ,$$
from which it follows that
$$\sup_{y\in Y}\inf_{x\in X}(f(x)+g(x,y))\leq \sup_Af+\epsilon<\inf_Bf\ .\eqno{(2.1)}$$
On the other hand, in view of $(c)$ and $(d)$, we have
$$\inf_Bf\leq \inf_{x\in B}(f(x)+\sup_{y\in Y}g(x,y))=\inf_{x\in B}\sup_{y\in Y}(f(x)+g(x,y))=
\inf_{x\in X}\sup_{y\in Y}(f(x)+g(x,y))\ .\eqno{(2.2)}$$
Now, the conclusion follows directly from $(2.1)$ and $(2.2)$.\hfill $\bigtriangleup$\par
\medskip
We also recall the following well-known fact:\par
\medskip
PROPOSITION 2.3. - {\it Let $A\subseteq {\bf R}^n$ be any open set and let $v\in L^1(A)\setminus
\{0\}$.\par
Then, one has
$$\sup_{\alpha\in C^{\infty}_0(A)}\int_{A}\alpha(x)v(x)dx=+\infty\ .$$}
\bigskip
{\bf 3. Proof and corollary of Theorem 1.1.}\par
\bigskip
 Let $\lambda\in ]a,b[$. 
The Lipschitz constant of the operator ${{\lambda}\over {\alpha_ {\Phi}L}}\Psi$ is equal
to ${{|\lambda|}\over {\alpha_ {\Phi}}}$ and so it is strictly less than ${{1}\over {\alpha_ {\Phi}}}$. Then, by Theorem 2.C, the multifunction
$y\to (\Phi+{{\lambda}\over {\alpha_ {\Phi}L}}\Psi)^{-1}(y)$  is lower semicontinuous (since it is Lipschitzian) and its values are non-empty and
connected (since they are absolute retracts). So, the operator $\Phi+\lambda\Psi$ is onto $Y$ and open.
Consequently, the operators $\Phi, \Psi$ satisfy the assumptions of Theorem 2.1, with $D=]a,b[$ and
$S_{\lambda}=Y$. Hence, the conclusion is directly ensured by Theorem 2.1.
 \hfill $\bigtriangleup$\par
\medskip
Let us notice explicitly the following corollary of Theorem 1.1.\par
\medskip
COROLLARY 3.1. - {\it Let $X, Y$ be two real Banach spaces, with $\hbox {\rm dim}(Y)\geq 2$, $\Phi:X\to Y$ a surjective continuous linear operator, $\Psi:X\to Y$ a non-constant Lipschitzian operator
with Lipschitz constant equal to $L$, $[a,b]$ a closed sub-interval of
$[-1,1]$.\par
Then, for each pair of continuous and convex functions $\theta:[0,+\infty[\to {\bf R}, \eta:[a,b]\to {\bf R}$, with $\theta$ strictly increasing, one has
$$\min\left \{\sup_{x\in X}\theta\left (\left \|\Phi(x)+{{a}\over {\alpha_ {\Phi}L}}\Psi(x)\right \|_Y\right )+\eta(a),
\sup_{x\in X}\theta\left (\left \|\Phi(x)+{{b}\over {\alpha_ {\Phi}L}}\Psi(x)\right \|_Y\right )+\eta(b)\right \}$$
$$=\sup_{x\in X}\inf_{\lambda\in [a,b]}\left (\theta\left (\left \|\Phi(x)+{{\lambda}\over {\alpha_ {\Phi}L}}\Psi(x)\right \|_Y\right )+\eta(\lambda)\right )\ .$$}
\smallskip
PROOF. Since dim$(Y)\geq 2$, the norm on $Y$ is a convex and sup-connected functional and hence so is $\theta(\|\cdot\|_Y)$.
Then, we can apply Theorem 2.1 taking
$\varphi(\cdot)=-\theta(\|\cdot\|_Y)$ and $\gamma=-\eta$, and the conclusion follows.\hfill $\bigtriangleup$
\medskip
REMARK 3.1. - Notice that Corollary 3.1 does not hold, in general, if $Y={\bf R}$. In this connection, it is enough
to take $X={\bf R}$, $\Phi(x)=x$, $\Psi(x)=|x|$, $[a,b]=[-1,1]$, $\theta(t)=t$, $\eta=0$. Hence, $L=\alpha_{\Phi}=1$ and we have
$$\sup_{x\in {\bf R}}\inf_{|\lambda|\leq 1}|x+\lambda |x||=0<+\infty=\inf_{|\lambda|\leq 1}\sup_{x\in {\bf R}}|x+\lambda |x||\ .$$
\bigskip
{\bf 4. Proof of Theorem 1.2.}\par
\bigskip
First, notice that, in view of the Rellich-Kondrachov theorem, for each $u\in W^{1,p}(\Omega)$, we have $\psi\circ u\in L^1(\Omega)$ and
the functional $u\to \int_{\Omega}\alpha(x)\psi(u(x))dx$ is sequentially weakly continuous. Moreover, by $(i)$, the functional $u\to \int_{\Omega}\varphi(x,u(x),\nabla u(x)dx$
 is sequentially weakly lower semicontinuous ([2], Theorem 4.6.8).
Now, let $V$ be a sequentially weakly closed subset of $W^{1,p}(\Omega)$ containing the constants and let $Y$ be a dense subset of $L^{\infty}(\Omega)$. Put
$$X=\left \{u\in V : \int_{\Omega}\varphi(x,u(x),\nabla u(x))dx<+\infty\right \}\ .$$
By $(ii)$, the constants belong to $X$.
Fix $\alpha\in L^{\infty}(\Omega)$. By $(i)$, there is $\delta>0$ such that
$$\varphi(x,\xi,\eta)-2\|\alpha\|_{L^{\infty}(\Omega)}|\psi(\xi)|\geq 0$$
for all $(x,\xi,\eta)\in \Omega\times{\bf R}\times{\bf R}^n$ with $|\xi|>\delta$. So, we have
$${{c}\over {2}}|\eta|^p-d-\|\alpha\|_{L^{\infty}(\Omega)}\sup_{|\xi|\leq\delta}|\psi(\xi)|\leq \varphi(x,\xi,\eta)+\alpha(x)\psi(\xi)$$
for all $(x,\xi,\eta)\in \Omega\times{\bf R}\times{\bf R}^n$ and, of course,
$$\lim_{|\xi\to +\infty}\inf_{(x,\eta)\in \Omega\times{\bf R}^n}(\varphi(x,\xi,\eta)+\alpha(x)\psi(\xi))=+\infty\ .$$
Consequently, in view of Proposition 2.1, we have, in $W^{1,p}(\Omega)$,
$$\lim_{\|u\|\to +\infty}\left (\int_{\Omega}\varphi(x,u(x),\nabla u(x))dx+\int_{\Omega}\alpha(x)\psi(u(x))dx\right )=+\infty \ .$$
This implies that, for each $r\in {\bf R}$, the set
$$\left \{u\in V: \int_{\Omega}\varphi(x,u(x),\nabla u(x))dx+\int_{\Omega}\alpha(x)\psi(u(x))dx\leq r\right \}$$
is weakly compact by reflexivity and Eberlein-Smulyan's theorem. Of course, we also have
$$\left \{u\in V: \int_{\Omega}\varphi(x,u(x),\nabla u(x))dx+\int_{\Omega}\alpha(x)\psi(u(x))dx\leq r\right \}$$
$$=\left \{u\in X: \int_{\Omega}\varphi(x,u(x),\nabla u(x))dx+\int_{\Omega}\alpha(x)\psi(u(x))dx\leq r\right \}\ .$$
Since the function $\int_{\Omega}\varphi(x,\cdot,0)dx$ is not quasi-convex,
 there are $\xi_1, \xi_2, \xi_3\in {\bf R}$, with $\xi_1<\xi_2<\xi_3$, such that
$$\max\left \{\int_{\Omega}\varphi(x,\xi_1,0)dx, \int_{\Omega}\varphi(x,\xi_3,0)dx\right \}<\int_{\Omega}\varphi(x,\xi_2,0)dx\ .$$
Now, observe that, if we put
$$A=\{\xi_1, \xi_3\}$$
and
$$B=\{\xi_2\}\ ,$$
and define $f:X\to {\bf R}$, $g:X\times Y\to {\bf R}$ by
$$f(u)=\int_{\Omega}\varphi(x,u(x),\nabla u(x))dx \ ,$$
$$g(u,\alpha)=\int_{\Omega}\alpha(x)(\psi(u(x))-\psi(\xi_2))dx$$
for all $u\in X$, $\alpha\in Y$, we clearly have
$$\sup_Af<\inf_Bf$$
and
$$\inf_{u\in B}\sup_{\alpha\in Y}g(u,\alpha)=0\ .$$
Since $\psi$ is strictly monotone, the numbers $\psi(\xi_1)-\psi(\xi_2)$ and $\psi(\xi_3)-\psi(\xi_2)$
have opposite signs. This clearly implies that
$$\sup_{\alpha\in Y}\inf_{u\in A}g(u,\alpha)\leq 0\ .$$
Furthermore, if $u\in X\setminus \{\xi_2\}$, again by strict monotonicity, $\psi\circ u\neq \psi(\xi_2)$, and so, since
$Y$ is dense in $L^{\infty}(\Omega)$, we have
$$\sup_{\alpha\in Y}g(u,\alpha)=+\infty\ .$$
Therefore, the sets $A, B$ and the functions $f, g$ satisfy the assumptions of Proposition 2.2. Consequently, we have
$$\sup_{\alpha\in Y}\inf_{u\in X}\left ( \int_{\Omega}\varphi(x,u(x),\nabla u(x))dx+\int_{\Omega}\alpha(x)\psi(u(x))dx\right )$$
$$<\inf_{u\in X}\sup_{\alpha\in Y}\left ( \int_{\Omega}\varphi(x,u(x),\nabla u(x))dx+\int_{\Omega}\alpha(x)\psi(u(x))dx\right )\ .$$
Now, the conclusion is a direct consequence of Theorem 2.B. When $Y=C^{\infty}_0(\Omega)$ the same proof as above holds in view of Proposition 2.3.\hfill $\bigtriangleup$\par
\medskip
We conclude presenting an application of Theorem 1.2 to the Neumann problem.\par
\medskip
 We denote by $\tilde{\cal A}$ the class of all Carath\'eodory functions
 $\psi:\Omega\times {\bf R}\to {\bf R}$ such that
$$\sup_{(x,\xi)\in\Omega\times {\bf R}}{{|\psi(x,\xi)|}\over
{1+|\xi|^{q}}}<+\infty\ ,$$
where  $0<q< {{pn-n+p}\over {n-p}}$ if $p<n$ and $0<q<+\infty$ if
$p=n$. While, when $n<p$, $\tilde{\cal A}$ stands for the class
of all Carath\'eodory functions $\psi:\Omega\times{\bf R}\to {\bf R}$. 
Given $\psi\in \tilde{\cal A}$, consider the following Neumann problem
$$\cases {-\hbox {\rm div}(|\nabla u|^{p-2}\nabla u)=
\psi(x,u)
 & in
$\Omega$\cr & \cr {{\partial u}\over {\partial \nu}}=0 & on
$\partial \Omega$\ ,\cr}\eqno{(P_{\psi})} $$
where $\nu$ is the outward unit normal to $\partial\Omega$.
 Let us recall
that a weak solution
of $(P_{\psi})$ is any $u\in W^{1,p}(\Omega)$ such that
 $$\int_{\Omega}|\nabla u(x)|^{p-2}\nabla u(x)\nabla v(x)dx
-\int_{\Omega}
\psi(x,u(x))v(x)dx=0$$
for all $v\in W^{1,p}(\Omega)$.\par
If $\psi\in\tilde{\cal A}$, we set $\Psi(x,\xi)=\int_0^{\xi}\psi(x,t)dt$. Clearly, $\Psi(x,\cdot)$ lies
in $\cal {A}$.\par
\medskip
THEOREM 4.1. - {\it Let $f, g:{\bf R}\to {\bf R}$ be two functions lying in $\tilde{\cal A}$ and satisfying the following
conditions:\par
\noindent
$(a_1)$\hskip 5pt the function $g$ has a constant sign and $int(g^{-1}(0))=\emptyset$\ ;\par
\noindent
$(a_2)$\hskip 5pt $\lim_{|\xi|\to +\infty}{{F(\xi)}\over {|G(\xi)|+1}}=+\infty$\ ;\par
\noindent
$(a_3)$\hskip 5pt the function $F-G$ is not quasi-convex\ .\par
Then, for each $\beta\in L^{\infty}(\Omega)$, with $\inf_{\Omega}\beta>0$, and for each convex set $Y\subseteq L^{\infty}(\Omega)$, dense
in $L^{\infty}(\Omega)$, there exists $\alpha\in Y$ such that the problem
$$\cases {-\hbox {\rm div}(|\nabla u|^{p-2}\nabla u)=\alpha(x)g(u)-\beta(x)f(u)
 & in
$\Omega$\cr & \cr {{\partial u}\over {\partial \nu}}=0 & on
$\partial \Omega$\ ,\cr}\eqno{(P)}$$
has at least three weak solutions.}\par
\smallskip
PROOF. Fix $\beta\in L^{\infty}(\Omega)$, with $\inf_{\Omega}\beta>0$, and a convex set $Y\subseteq L^{\infty}(\Omega)$, dense
in $L^{\infty}(\Omega)$. We are going to apply Theorem 1.2, defining $\varphi, \psi$  by
$$\varphi(x,\xi,\eta)={{1}\over {p}}|\eta|^p+\beta(x)(F(\xi)-G(\xi))$$
and
$$\psi(\xi)=-G(\xi)$$
for all $(x,\xi,\eta)\in \Omega\times {\bf R}\times {\bf R}^n$. It is immediate to realize that, by $(a_1)-(a_3)$, the above $\varphi, \psi$
satisfy the assumptions of Theorem 1.2. Of course, the set $Y-\beta$ is convex and dense in $L^{\infty}(\Omega)$. Then, Theorem 1.2
ensures the existence of $\alpha\in Y$ such that the functional
$$u\to \int_{\Omega}\varphi(x,u(x),\nabla u(x))dx+\int_{\Omega}(\alpha(x)-\beta(x))\psi(u(x))dx$$
$$={{1}\over {p}}\int_{\Omega}|\nabla u(x)|^pdx+\int_{\Omega}\beta(x)F(u(x))dx-\int_{\Omega}\alpha(x)G(u(x))dx$$
has at least two global minima in $W^{1,p}(\Omega)$. But, by classical results, such a functional is $C^1$ and satisfies the Palais-Smale condition
and hence, by Corollary 1 of [3], has at least three critical points which are weak solutions of problem $(P)$.\hfill $\bigtriangleup$
\medskip
In conclusion, we want to remark a feature of Theorem 2.A, the main tool that we used to prove Theorem 1.1: the second variable of $f$ runs over a real interval. An important contribution to this kind of results has been provided by H. Tuy in [17].

\bigskip
\bigskip
{\bf Acknowledgement.} The author has been supported by the Gruppo Nazionale per l'Analisi Matematica, la Probabilit\`a e 
le loro Applicazioni (GNAMPA) of the Istituto Nazionale di Alta Matematica (INdAM) and by the Universit\`a degli Studi di Catania, ``Piano della Ricerca 2016/2018 Linea di intervento 2". \par
\bigskip
\bigskip
\bigskip
\bigskip
\vfill\eject
\centerline {\bf References}\par
\bigskip
\bigskip
\noindent
[1]\hskip 5pt DENKOWSKI, S. MIG\'ORSKI and N. S. PAPAGEORGIOU,
{\it An Introduction to Nonlinear Analysis: Theory}, Kluwer Academic
Publishers, 2003.\par
\smallskip
\noindent
[2]\hskip 5pt Z. DENKOWSKI, S. MIG\'ORSKI and N. S. PAPAGEORGIOU,
{\it An Introduction to Nonlinear Analysis: Applications}, Kluwer Academic
Publishers, 2003.\par
\smallskip
\noindent
[3]\hskip 5pt  P. PUCCI and J. SERRIN, {\it A mountain pass theorem},
J. Differential Equations, {\bf 60} (1985), 142-149.\par
\smallskip
\noindent
[4]\hskip 5pt B. RICCERI, {\it Structure, approximation et d\'ependance
continue des solutions de certaines \'equations non lin\'eaires},
C. R. Acad. Sci. Paris, S\'erie I, {\bf 305} (1987), 45-47.\par
\smallskip
\noindent
[5]\hskip 5pt B. RICCERI, {\it Some topological mini-max theorems via
an alternative principle for multifunctions}, Arch. Math. (Basel),
{\bf 60} (1993), 367-377.\par
\smallskip
\noindent
[6]\hskip 5pt B. RICCERI, {\it On a topological minimax theorem and
its applications}, in ``Minimax theory and applications'', B. Ricceri
and S. Simons eds., 191-216, Kluwer Academic Publishers, 1998.\par
\smallskip
\noindent
[7]\hskip 5pt B. RICCERI, {\it A further improvement of a minimax theorem of
Borenshtein and Shul'man}, J. Nonlinear Convex Anal., {\bf 2} (2001),
279-283.\par
\smallskip
\noindent
[8]\hskip 5pt B. RICCERI, {\it Minimax theorems for limits of 
parametrized functions
having at most one local minimum lying in a certain set}, Topology Appl., 
{\bf 153} (2006), 3308-3312.\par
\smallskip
\noindent
[9]\hskip 5pt B. RICCERI, {\it Recent advances in minimax theory and
applications}, in ``Pareto Optimality, Game Theory and Equilibria'', 
A. Chinchuluun, P.M. Pardalos, A. Migdalas, L. Pitsoulis eds., 23-52,
Springer, 2008.\par
\smallskip
\noindent
[10]\hskip 5pt B. RICCERI, 
{\it Nonlinear eigenvalue problems},  
in ``Handbook of Nonconvex Analysis and Applications'' 
D. Y. Gao and D. Motreanu eds., 543-595, International Press, 2010.\par
\smallskip
\noindent
[11]\hskip 5pt B. RICCERI, {\it A strict minimax inequality criterion and some of its consequences}, Positivity, {\bf 16} (2012), 455-470.\par
\smallskip
\noindent
[12]\hskip 5pt B. RICCERI, {\it Energy functionals of Kirchhoff-type problems having multiple global minima}, Nonlinear Anal., {\bf 115} (2015),  
130-136.\par
\smallskip
\noindent
[13]\hskip 5pt B. RICCERI, {\it A minimax theorem in infinite-dimensional topological vector spaces}, Linear Nonlinear Anal., {\bf 2}
(2016), 47-52.\par
\smallskip
\noindent
[14]\hskip 5pt B. RICCERI, {\it Miscellaneous applications of certain minimax theorems I}, Proc. Dynam. Systems Appl., {\bf 7} (2016),
198-202.\par
\smallskip
\noindent
[15]\hskip 5pt B. RICCERI, {\it On a minimax theorem: an improvement, a new proof and an overview of its applications},
Minimax Theory Appl., {\bf 2} (2017), 99-152.\par
\smallskip
\noindent
[16]\hskip 5pt R. T. ROCKAFELLAR, {\it  Integral functionals, normal integrands and measurable selections}, Lecture Notes in Math., Vol. 543,  157-207,
Springer, Berlin, 1976.\par
\smallskip
\noindent
[17]\hskip 5pt H. TUY, {\it A new topological minimax theorem with application},
J. Global Optim., {\bf 50} (2011), 371-378. \par
\bigskip
\bigskip
\bigskip
\bigskip
Department of Mathematics and Informatics\par
University of Catania\par
Viale A. Doria 6\par
95125 Catania, Italy\par
{\it e-mail address}: ricceri@dmi.unict.it

\bye